\documentclass[a4paper]{article}

\usepackage{arxiv}
\usepackage[utf8]{inputenc} 
\usepackage[T1]{fontenc}    
\usepackage{hyperref}       
\usepackage{url}            
\usepackage{booktabs}       
\usepackage{amsfonts}       
\usepackage{nicefrac}       
\usepackage{microtype}      
\usepackage{lipsum}		
\usepackage{graphicx}
\usepackage{natbib}
\usepackage{doi}

\usepackage{amsmath}
\usepackage{amsfonts}
\usepackage{amssymb}
\usepackage{amsthm}
\usepackage{newlfont}
\usepackage{amsbsy}
\usepackage{bm}
\usepackage[all]{xy}
\usepackage{textcomp}
\usepackage{color}
\usepackage{epsfig}
\usepackage{dsfont}
\usepackage{latexsym}
\usepackage{array}
\usepackage{longtable}
\usepackage{enumerate}
\usepackage{multicol}
\usepackage{mathrsfs}
\usepackage{cancel}
\usepackage{wasysym}

\newtheorem{teor}{Theorem}
\newtheorem{lemma}{Lemma}
\newtheorem{coro}{Corollary}
\newtheorem{prop}{Proposition}

\theoremstyle{definition}

\newtheorem*{ejem*}{Examples}
\newtheorem{defin}{Definition}
\newtheorem*{demos}{Proof}
\newtheorem{remark}{Remark}
\theoremstyle{remark}

\newcommand{\complex}{\mathbf{\mathbb{C}}}
\newcommand{\reales}{\mathbf{\mathbb{R}}}
\newcommand{\racio}{\mathbf{\mathbb{Q}}}

\newcommand{\natu}{\mathbf{\mathbb{N}}}
\newcommand{\llave}[1]{\left\{ #1\right\}}

\newcommand{\Sum}[2]{\sum\limits_{#1}^{#2}}

\newcommand{\corch}[1]{\left[ #1\right]}
\newcommand{\paren}[1]{\left( #1\right)}
\newcommand{\QED}{\hfill \ensuremath{\Box}}

\title{Characterization of Rational Solutions of a KdV-Like Equation}
\date{\today}
\author{ {Brian D. Vasquez} \\
	Khalifa University, Abu Dhabi\\
	\texttt{brian.campos@ku.ac.ae} \\
}

\begin{document}

\maketitle

\begin{abstract}
We characterize the rational solutions to a KdV-like equation which are generated from polynomial solutions to the corresponding generalized bilinear equation. 
We use a particular class of polynomials satisfying a quadratic difference equation to obtain that there are no solutions of the bilinear equation with degree (in the spatial variable) greater than $5$.
As a byproduct, we answer positively  a  conjecture of Yi Zhang and Wen-Xiu Ma about these solutions.\\



\end{abstract}

\section*{Introduction}

Rational solutions of nonlinear partial differential equations have been considered as a simple way to understand their associated flows. Two remarkable cases are the KdV and the KP equations. In the KdV case, it is known after the work of Airault, McKean, and Moser \cite{AMM} that the rational solutions $u(x,t)$ which decay to zero are of the form
\begin{equation*}
    u(x,t)=\frac{1}{2}\Sum{i=1}{n}\frac{1}{(x-x_{i}(t))^2}.
\end{equation*}
On the other hand, the method developed by Krichever \cite{kricheverrat} for constructing rational solutions of the Zakharov-Shabat equations allows to understand the motion of the poles $x_{1}(t), ... , x_{n}(t)$ of rational solutions of the KP equation. This method uses the Calogero-Moser theory for a system of $n$ particles on a line with the Hamiltonian
\begin{equation*}
    H(p,x)=\frac{1}{2}\Sum{i=1}{n}p_{i}^2+\Sum{i<j}{}2(x_{i}-x_{j})^{-2}.
\end{equation*}

Furthermore, rational solutions of nonlinear partial differential equations yield bispectral AKNS/ZS operators \cite{ZUBELLI199271,MR1857803}. 
In fact, if $q(x,t_{2},...,t_{m})$ and $r(x,t_{2},...,t_{m})$ are certain rational solutions of the AKNS hierarchy then the matrix differential operator
\begin{equation*}
    L=	\paren{\begin{matrix} 
\partial_{x} & -q \\
	r & -\partial_{x}
	\end{matrix}},
\end{equation*}
is bispectral, i.e., there exist eigenfunctions $\psi(x,z)$ of $L$ such that $B(z,\partial_{z})\psi=\theta(x)\psi$ for some matrix differential operator $B=B(z,\partial_{z})$ which is independent of the variable $x$ and $\theta$ is a nonconstant function of $x$.
In \cite{Zubelli2000} rational solutions which decay at infinity and are bispectral potentials for the Schr\"{o}dinger operator were characterized. Those functions are known as master symmetries for KdV.

As an application of the  rational solutions of nonlinear partial differential equations, we have the rogue wave solutions as a model to be aware of the unexpected nonlinear wave phenomena in the ocean called rogue waves \cite{kharif2008rogue}, \cite{ocean}.

We consider the generalized bilinear differential equation of KdV type
\begin{equation}\label{bildiff}
    2f_{xt}f-2f_{t}f_{x}+6(f_{xx})^2=0.
\end{equation}

By taking the logarithmic derivative transformation 
\begin{equation}\label{logtrans}
    u=2\partial_{x} log (f),
\end{equation}
we obtain the KdV-like equation 
\begin{equation}\label{KdVlike}
    u_{t}+\frac{3}{2}(u_{x})^2+\frac{3}{2}u^2u_{x}+\frac{3}{8}u^4=0.
\end{equation}

Note that a polynomial $f$ solves Equation \eqref{bildiff} if and only if  $u$ is a rational solution of the KdV-like Equation \eqref{KdVlike}.
In  \cite{bilinearMa} generalizations of the bilinear differential operator of Hirota form (see \cite{zbMATH02117215}) were first studied.

To define the generalized operators we fix two integers $p,k \in \natu$ and consider
\begin{equation*}
    D_{p,x_{1}}^{n_{1}}\cdot \cdot \cdot D_{p,x_{k}}^{n_{k}}f\cdot g=
    \prod_{i=1}^{k} \paren{\frac{\partial}{\partial x_{i}}+\alpha_{p}\frac{\partial}{\partial x_{i}^{'}}}^{n_{i}}f(x_{1},...,x_{k})g(x_{1}^{'},...,x_{k}^{'})
    \mid_{x_{1}^{'}=x_{1}, ... ,x_{k}^{'}=x_{k} },
\end{equation*}
where $n_{1}, ... , n_{k}$ are arbitrary nonnegative integers and the powers of $\alpha_{p}$ are defined for an integer $m$ as
\begin{equation*}
    \alpha_{p}^{m}=(-1)^{r_{p}(m)}, \hspace{0.1cm} \text{if} \hspace{0.1cm} m\equiv r_{p}(m) \hspace{0.1cm} (mod \hspace{0.1cm} p) \hspace{0.1cm} \text{with} \hspace{0.1cm} 0\leq r_{p}(m)<p.
\end{equation*}

In particular, if we consider the case $p=3$, $k=2$, then $$D_{3,x}D_{3,t}f \cdot f=2f_{xt}f-2f_{t}f_{x}, \hspace{0.1cm}  D_{3,x}^4
f\cdot f=6(f_{xx})^2.$$ Therefore, $$\paren{D_{3,x}D_{3,t}+D_{3,x}^4}f \cdot f=2f_{xt}f-2f_{t}f_{x}+6(f_{xx})^2,$$ which is exactly the operator defining the Equation \eqref{bildiff}. 

Furthermore, if we consider the Hirota case, i.e.,  $p=k=2$, we obtain
$$D_{2,x}D_{2,t} f\cdot f=2f_{xt}f-2f_{t}f_{x}, \hspace{0.1cm}
D_{2,x}^4 f\cdot f=2f_{xxxx}f-8f_{xxx}f_{x}+ 6(f_{xx})^2.$$
Therefore, $$\paren{D_{2,x}D_{2,t}+D_{2,x}^4}f \cdot f=2f_{xt}f-2f_{t}f_{x}+2f_{xxxx}f-8f_{xxx}f_{x}+ 6(f_{xx})^2,$$
which is exactly the bilinear equation associated to the KdV equation. See \cite{zbMATH02117215}.

In the present work we characterize the polynomial solutions of Equation \eqref{bildiff}. 
We show that the condition $\deg_{x}(f) \leq 4$ is sufficient to obtain polynomial solutions of Equation \eqref{bildiff} by giving their explicit expressions. 

In order to characterize the algebro-geometric structure of the set of solutions we start by defining a classical notation used in algebraic geometry.

\begin{defin}
Let $X$ be an arbitrary set, $\mathbb{K}$ a field and $f:X \rightarrow \mathbb{K}$ an arbitrary function. We define 
\begin{equation*}
    V(f)=\llave{x\in X\mid f(x)=0}
\end{equation*}
the zero set of $f$.
\end{defin}

\begin{defin}
Let $U\subset \reales^d$ open, $f_{1},...,f_{n}\in C^{n-1}(U)$. For $1\leq k \leq d$ we define the function
$W_{x_{k}}(f_{1},...,f_{n})\in C(U)$ 
\begin{equation*}
    W_{x_{k}}(f_{1},...,f_{n})=\det \corch{\paren{\frac{\partial^{j-1} f_{i}}{\partial x_{k}^{j-1}}}_{i,j\in \llave{1,...,n}}}.
\end{equation*}
We call $W_{x_{k}}$ the partial wronskian with respect to $x_{k}$.
\end{defin}

\begin{defin}
Let $A$ be a graded ring. We will denote the homogeneous component of degree $k\geq 0$ by $A_{k}$. 
\end{defin}

We consider the nonlinear operator $T:\complex[x,t]\xrightarrow[]{} \complex[x,t]$, $Tf(x,t)=W_{t}(f_{x},f)+3(f_{xx})^2$. Note that the polynomial solutions of the Equation \eqref{bildiff} are given by the set $V(T)$.

The plan of this article is as follows: In Section \ref{sec1}, we fix $m\in \natu$ and decompose the operator $T$ as a sum of two operators $B_{(m)},R_{(m)}:\complex[x,t]\xrightarrow[]{} \complex[x,t]$ such that $T=B_{(m)}+R_{(m)}$. In Section \ref{sec2}, we study an extension in the variable $m$ of the time-depending polynomials $P_{k}$, $0 \leq k \leq m$ defining the fundamental solutions of the operator $B_{(m)}$ and its leading polynomials in $m$. By fundamental solution we mean a monic polynomial $f\in V(B_{(m)})$ which generates $V(B_{(m)})$ as $\complex[t]-$module. In Section \ref{sec3}, we consider the operator remainder $R_{(m)}$ and  time-depending polynomials $R_{j}$, $m+1\leq j\leq 2m-1$ which define it. In Subsection \ref{sub3.1}, we find a general formula for the leading function of $R_{j}$ as a polynomial in the leading polynomials of the time-depending polynomials $P_{k}$ defining the fundamental solutions of the operator $B_{(m)}$.  In Subsection \ref{sub3.2}, we consider the algebraic sets defined through the leading functions of some remarkable time-depending remainder polynomials and show that their irreducible components are hyperplanes. Finally, in Section \ref{sec4}, we use those algebraic sets to show that we do not have polynomial solutions of the differential equation \eqref{bildiff} with $\deg_{x}(f)>5$. In particular, we obtain two classes of fundamental solutions of the bilinear differential equation \eqref{bildiff}. By taking the transformation in Equation \eqref{logtrans} we obtain two classes of rational solutions to the KdV-like Equation \eqref{KdVlike}.

This give a positive answer to the conjecture proposed in \cite{rationalzhangMa}.

\section{Decomposition of the Quadratic Operator}\label{sec1}

In this section, we face the problem of understanding the operator $T$ through the set $V(T)$. We start with a lemma about the decomposition of $T$.

\begin{lemma}
For every $m\in \natu$, there exist unique operators $B_{(m)},R_{(m)}:\complex[x,t]\xrightarrow[]{} \complex[x,t]$ such that for every $f\in \complex[x,t]$ we have $\deg_{x}(R_{(m)}f)\leq m-2$, $B_{(m)}f\in x^{m-1}\complex[x,t]$ and $T=B_{(m)}+R_{(m)}$. 
\end{lemma}
\begin{demos}
 Note that we have a natural decomposition $T=B_{(m)}+R_{(m)}$ with 
\begin{equation*}
    B_{(m)}f(x,t)=\Sum{p=m-1}{2n-1}\frac{\partial^{p}Tf}{\partial x^p}(0,t)x^p \hspace{0.1cm} \text{and} \hspace{0.1cm}
    R_{(m)}f(x,t)=\Sum{p=0}{m-2}\frac{\partial^{p}Tf}{\partial x^p}(0,t)x^p,
\end{equation*}
for every polynomial $f\in \complex[x,t]$, $\deg_{x}(f)\leq n$.\QED
\end{demos}

\begin{teor}\label{RandB}
Let $m\in \natu$. If $\deg_{x}(f)\leq m$, then
\begin{equation}\label{remainder}
    R_{(m)}f(x,t)=\Sum{p=0}{m-2}\llave{\Sum{j=0}{m}jW(f_{j},f_{p+1-j})(t)
    +3\Sum{j=0}{m}j(j-1)(p+4-j)(p+3-j)f_{j}(t)f_{p+4-j}(t)}x^p
\end{equation}
and 
\begin{equation}\label{bil}
    B_{(m)}f(x,t)=\Sum{p=2m-3}{2m-1}\Sum{j=0}{m}jW(f_{j},f_{p+1-j})(t)x^p
    \end{equation}
    \begin{equation*}
    +\Sum{p=m-1}{2m-4}\llave{\Sum{j=0}{m}jW(f_{j},f_{p+1-j})(t)
    +3\Sum{j=0}{m}j(j-1)(p+4-j)(p+3-j)f_{j}(t)f_{p+4-j}(t)}x^p,
\end{equation*}
for each $f\in \complex[x,t]$.
\end{teor}
\begin{demos}
If $f(x,t)=\Sum{j=0}{m}f_{j}(t)x^{j}$ then $f_{t}(x,t)=\Sum{j=0}{m}f_{j}^{'}(t)x^{j}$,  $f_{x}(x,t)=\Sum{j=0}{m}jf_{j}(t)x^{j-1}$ and 
$f_{xt}(x,t)=\Sum{j=0}{m}jf_{j}^{'}(t)x^{j-1}$, $f_{xx}(x,t)=\Sum{j=0}{m}j(j-1)f_{j}(t)x^{j-2}$ we have the expression
\begin{equation*}
    Tf(x,t)=W_{t}(f_{x},f)+3(f_{xx})^2
    =\Sum{p=2m-3}{2m-1}\Sum{j=0}{m}jW(f_{j},f_{p+1-j})(t)x^p
\end{equation*} 
\begin{equation*}
    +\Sum{p=0}{2m-4}\llave{\Sum{j=0}{m}jW(f_{j},f_{p+1-j})(t)
    +3\Sum{j=0}{m}j(j-1)(p+4-j)(p+3-j)f_{j}(t)f_{p+4-j}(t)}x^p.
\end{equation*}
Thus the operators $B_{(m)}$ and $R_{(m)}$ are given by Equations \eqref{remainder} and \eqref{bil}. \QED
\end{demos}

\begin{coro}
For every $m\in \natu$ we have  $V(T)=V(B_{(m)})\cap V(R_{(m)})$.
\end{coro}

We now fix $m\in \natu$ and denote $B_{(m)}$ and $R_{(m)}$ simply by $B$ and $R$. Note that since $V(T)=V(B)\cap V(R)$ we can understand $V(T)$ by means of $V(B)$. 


\begin{prop}\label{charac}
Let $f\in V(B)$, $f(x,t)=\Sum{j=0}{m}f_{j}(t)x^{j}$. Then for every $0\leq k \leq m$ there exists $P_{m,k}\in \complex[t]$ such that $f_{k}(t)=P_{m,k}(t)f_{m}(t)$ and 
\begin{equation*}
    P_{m,k}^{'}(t)=\frac{1}{m-k}\left\{\Sum{j=k+1}{m-1}jW(P_{m,k+m-j},P_{m,j})(t)\right.
  \end{equation*}
\begin{equation}\label{polyPK}
     \left. +3\Sum{j=k+3}{m}j(j-1)(k+m+3-j)(k+m+2-j)P_{m,j}(t)P_{m,k+m+3-j}(t) \right\}.
\end{equation}
\end{prop}
\begin{demos}
Let us assume that $f\in \complex[x,t]$, $f(x,t)=\Sum{j=0}{m}f_{j}(t)x^{j}$ satisfies $Bf=0$ and look for the restrictions over $f$.
In fact, by Theorem \ref{RandB} we have  $Bf=0$ if and only if 
\begin{equation*}
   \Sum{j=0}{m} jW(f_{j},f_{p+1-j})(t)=0 \hspace{0.1cm} \text{for every} \hspace{0.1cm} 2m-3\leq p \leq 2m-1,
\end{equation*}
\begin{equation*}
\Sum{j=p+1-m}{m}jW(f_{j},f_{p+1-j})(t)
    +3\Sum{j=p+4-m}{m}j(j-1)(p+4-j)(p+3-j)f_{j}(t)f_{p+4-j}(t)=0 \hspace{0.1cm} \text{for every} \hspace{0.1cm} m-1\leq p \leq 2m-4.
\end{equation*}

Therefore, 
\begin{equation*}
    W(f_{p+1-m},f_{m})(t)=\frac{1}{p+1-2m}
    \left\{\Sum{j=p+2-m}{m-1}jW(f_{j},f_{p+1-j})(t)\right.
    \end{equation*}
\begin{equation}\label{polyfk}
    \left. -3\Sum{j=p+4-m}{m}j(j-1)(p+4-j)(p+3-j)f_{j}(t)f_{p+4-j}(t) \right\} \hspace{0.1cm} \text{for every} \hspace{0.1cm} m-1\leq p \leq 2m-4.
\end{equation}
We claim that $f_{j}(t)=P_{m,j}(t)f_{m}(t)$ for some polynomial $P_{m,j}(t)\in \complex[t]$. We write $m-j=3q+r$ for $r\in \llave{0,1,2}$ and we make induction over $q$ for $r\in \llave{0,1,2}$ arbitrary. 
Note that 
\begin{itemize}
    \item If $p=2m-1$ we obtain an identity $W(f_{j},f_{j})(t)=0$.
    \item If $p=2m-2$ then $\Sum{j=0}{m} jW(f_{j},f_{2m-1-j})(t)=\Sum{j=m-1}{m} jW(f_{j},f_{2m-1-j})(t)\\
    =-W(f_{m-1},f_{m})(t)=0$ and we have 
    $f_{m-1}(t)=c_{m-1}f_{m}(t)$ for some $c_{m-1}\in \complex$.
    
    \item If $p=2m-3$ then $\Sum{j=0}{m} jW(f_{j},f_{2m-2-j})(t)=\Sum{j=m-2}{m} jW(f_{j},f_{2m-2-j})(t)\\
    =-2W(f_{m-2},f_{m})(t)=0$ and we have 
    $f_{m-2}(t)=c_{m-2}f_{m}(t)$ for some $c_{m-2}\in \complex$.
\end{itemize}
In particular, the case $q=0$ is okay since $P_{m,j}$ is constant for $j=m-2,m-1,m$.

Assume this assertion for $0\leq l \leq q$, $3(q+1)\leq m$. Then we have $f_{j}(t)=P_{m,j}(t)f_{m}(t)$ for $m-3q-2\leq j\leq m$. Now we consider the case 
$q+1$. 

If we put $p=2m-3q-4-r$ in \eqref{polyfk} for $r\in \llave{0,1,2}$ then 
\begin{equation*}
    \paren{\frac{f_{m-3q-3-r}}{f_{m}}}^{'}(t)=\frac{1}{-3q-3-r} \frac{1}{f_{m}(t)^2}\left\{\Sum{j=m-3q-2-r}{m-1}jW(f_{j},f_{2m-3q-3-r-j})(t)\right.
    \end{equation*}
\begin{equation*}    
   \left. -3 \Sum{j=m-3q-r}{m}j(j-1)(2m-3q-r-j)(2m-3q-r-j-1)f_{j}(t)f_{2m-3q-r-j}(t)\right\}
\end{equation*}
\begin{equation*}
    =\frac{1}{-3q-3-r} \left\{\Sum{j=m-3q-2-r}{m-1}jW(P_{m,j},P_{m,2m-3q-3-r-j})(t)\right.
\end{equation*}
\begin{equation*}    
   \left. -3 \Sum{j=m-3q-r}{m}j(j-1)(2m-3q-r-j)(2m-3q-r-j-1)P_{m,j}(t)P_{m,2m-3q-r-j}(t)\right\} .
\end{equation*}
Thus, $f_{m-3q-3-r}(t)=P_{m,m-3q-3-r}(t)f_{m}(t)$ with 
\begin{equation*}
  P_{m,m-3q-3-r}^{'}(t)  =\frac{1}{3(q+1)+r} \left\{\Sum{j=m-3q-2-r}{m-1}jW(P_{m,2m-3q-3-r-j},P_{m,j})(t)\right.
\end{equation*}
\begin{equation*}    
   \left. +3 \Sum{j=m-3q-r}{m}j(j-1)(2m-3q-r-j)(2m-3q-r-j-1)P_{m,j}(t)P_{m,2m-3q-r-j}(t)\right\} .
\end{equation*}
This is exactly the case $q+1$ for $r\in \llave{0,1,2}$ arbitrary and we obtain the formula \eqref{polyPK}.
Then the proof of the assertion follows by induction. \QED
\end{demos}

A simple consequence of this theorem is that we can consider monic polynomial solutions of \eqref{bildiff} in the variable $x$.

At first glance the set $V(B)$ wouldn't be a $\complex$-vector space since $B$ is not linear. However, it has a nice algebraic structure as we shall see in the following result.

\begin{teor}
 $V(B_{(m)})\cap \bigoplus_{k=0}^{m} (\complex[t])[x]_{k}$ is a cyclic $\complex[t]$-module generated by the monic polynomial $\overline{f}=\Sum{j=0}{m}P_{m,j}(t)x^j$. 
\end{teor}
\begin{demos}
Let $\overline{f}(x,t)=\Sum{j=0}{m}P_{m,j}(t)x^j$ with $\llave{P_{m,j}}_{0\leq j\leq m}$ satisfying the system of equations \eqref{polyPK}. Then the Proposition \ref{charac} implies that $V(B)\cap \bigoplus_{k=0}^{m} (\complex[t])[x]_{k}=\complex[t]\cdot \overline{f}$, i.e., $V(B)\cap \bigoplus_{k=0}^{m} (\complex[t])[x]_{k}$ is a cyclic $\complex[t]$-module generated by $\overline{f}$.
\end{demos}\QED

From now on, we will call the monic polynomial $\overline{f}$ the fundamental solution of degree $m$. \\

As a corollary we have 

\begin{coro}
$V(B)\cap \bigoplus_{k=0}^{m} (\complex[t])[x]_{k}$ is a $\complex$-vector space.
\end{coro}

Now we have an upper bound for the degree of the polynomial $P_{m,k}$, $0\leq k\leq m$.
\begin{lemma}\label{lemmadeg}
If $m-k=3q+r$ for $r\in \llave{0,1,2}$ then $\deg\paren{P_{m,k}}\leq q$.
\end{lemma}
\begin{demos}
The proof is by induction over $q$ for $r\in \llave{0,1,2}$ arbitrary.
For $q=0$ we have $k=m-r$, $r\in \llave{0,1,2}$. Since $P_{m,k}$ is constant for $k\in \llave{m-1,m-2,m-3}$ we have that $\deg\paren{P_{m,k}}=0=\frac{m-k-r}{3}$.

Assume that $\deg\paren{P_{m,j}}\leq l$, $m-j=3l+r$, $r\in \llave{0,1,2}$, $0\leq l\leq q-1$ and note that by Theorem \ref{charac} we have 

\begin{equation}\label{deg}
\deg\paren{P_{m,k}^{'}}\leq 
\max\left\{\max_{k+1\leq j \leq m-1}\llave{\deg\paren{P_{m,k+m-j}},\deg\paren{P_{m,j}}}-1,\right.
\end{equation}
\begin{equation*}
    \left. \max_{k+3\leq j\leq m}\llave{\deg\paren{P_{m,j}},\deg\paren{P_{m,k+m-j}}}\right\}.
\end{equation*}
Assume that $m-k=3q+r$. If $m-j=3n+s$, $s\in \llave{0,1,2}$ then $j-k=(m-k)-(m-j)=3(q-n)+r-s$. Since  $r,s\in \llave{0,1,2}$ we have 
 $r-s\in \llave{-2,-1,0,1,2}$ and $r-s=-3a+b$, $a\in \llave{0,1}$, $b\in \llave{0,1,2}$, therefore $j-k=3(q-n-a)+b$.

Using the induction hypothesis we obtain $\deg\paren{P_{m,j}}\leq n$ and $\deg\paren{P_{m,k+m-j}}\leq q-n-a$, therefore
\begin{equation*}
    \deg\paren{P_{m,k+m-j}}+\deg\paren{P_{m,j}}-1\leq q-1-a, \hspace{0.1cm} \text{for every} \hspace{0.1cm} k+1\leq j\leq m-1.
\end{equation*}

On the other hand, $\deg\paren{P_{m,k+m+3-j}}\leq q-n-1-a$ since $j-k-3=3(q-n-1-a)+b$. Hence, 
\begin{equation*}
    \deg\paren{P_{m,j}}+ \deg\paren{P_{m,k+m+3-j}}\leq q-1-a,\hspace{0.1cm} \text{for every} \hspace{0.1cm} k+3\leq j\leq m-1.
\end{equation*}
Thus, the Equation \eqref{deg} implies that 
\begin{equation*}
    \deg\paren{P_{m,k}^{'}}\leq q-1-a, \hspace{0.1cm} \text{in other words} \hspace{0.1cm} \deg\paren{P_{m,k}}\leq q-a\leq q.
\end{equation*}
The proof of the assertion follows by induction. \QED
\end{demos}

\section{The Polynomial Extension in the Index Variable and its Leading Coefficient }\label{sec2}

We now consider the function $P_{m,m-k}(t)$ and note that it has a polynomial extension in the variable $m$ as we shall see in the following 
result.
\begin{prop}\label{recurQ}
If $Q_{k}(t)=P_{m,m-k}(t)$ then 
\begin{equation*}
    Q_{k}(t)=\frac{1}{k}\left\{\Sum{l=1}{k-1}(m-l)\int W_{t}(Q_{k-l}(t),Q_{l}(t))dt\right.
\end{equation*}   
\begin{equation*}
    \left.+3\Sum{l=0}{k-3}(m-l)(m-l-1)(m+l-k+3)(m+l-k+2)\int Q_{l}(t)Q_{k-l-3}(t)dt\right\}+c_{k}
\end{equation*}
for an arbitrary $c_{k}\in \complex.$
\end{prop}
\begin{demos}
By Theorem \ref{charac} we have that 
\begin{equation*}
    \frac{d Q_{k}}{dt}(t)=
    \frac{1}{k}\left\{\Sum{j=m-k+1}{m-1}jW(P_{m,2m-k-j},P_{m,j})(t)\right.
  \end{equation*}
\begin{equation*}
     \left. +3\Sum{j=m-k+3}{m}j(j-1)(2m-k+3-j)(2m-k+2-j)P_{m,j}(t)P_{m,2m-k+3-j}(t) \right\}
\end{equation*}
\begin{equation*}
=\frac{1}{k}\left\{\Sum{l=1}{k-1}(m-l)W(P_{m,m+l-k},P_{m,m-l})(t)\right.
  \end{equation*}
\begin{equation*}
     \left. +3\Sum{l=0}{k-3}(m-l)(m-l-1)(m+l-k+3)(m+l-k+2)P_{m,m-l}(t)P_{m,m+l-k+3}(t) \right\}
\end{equation*}
\begin{equation*}
=\frac{1}{k}\left\{\Sum{l=1}{k-1}(m-l)W_{t}(Q_{k-l}(t),Q_{l}(t))\right.
  \end{equation*}
\begin{equation*}
     \left. +3\Sum{l=0}{k-3}(m-l)(m-l-1)(m+l-k+3)(m+l-k+2)Q_{l}(t)Q_{k-l-3}(t) \right\}
\end{equation*}
Taking the integral with respect to $t$ we conclude the proof.\QED
\end{demos}

This proposition tells us that $Q_{k}\in \complex[m,t,c_{1},...,c_{k}]$ for every $k\in \natu.$ 
Note that by the Lemma \ref{lemmadeg} we have that $\deg_{t} (Q_{k})\leq q$ with $k=3q+r$, $r\in \llave{0,1,2}.$ However, if we consider $Q_{k}$ as an element in $\complex[m,t,c_{1},...,c_{k}]$ we obtain that in fact $\deg_{t} (Q_{k})= q$. We can write  $Q_{k}(m,t,c_{1},...,c_{k})=\Sum{j=0}{q}x_{k,j}(m,c_{1},...,c_{k})t^{j}$ and look for conditions on the leading polynomial $x_{k,q}$ as we give in the following theorem.

\begin{teor}
The leading polynomial of $Q_{k}$ satisfies 
\begin{enumerate}
    \item $x_{k,q}\in \racio[m]$ with leading coefficient $\frac{1}{q!}$, if $k\equiv 0$ (mod 3), \label{1}
    \item $x_{k,q}\in \racio[m]c_{1}$ with leading coefficient $\frac{c_{1}}{q!}$, if $k\equiv 1$ (mod 3),\label{2}
    \item $x_{k,q}\in \racio[m]c_{1}^2+\racio[m]c_{2}$ with leading coefficient $\frac{c_{2}}{q!}$, if $k\equiv 2$ (mod 3).\label{3}
    \end{enumerate}
    Furthermore, we have $x_{0,r}=c_{r}$,  for $r\in \llave{0,1,2}$ with the recursion 

    \begin{equation*}
x_{k,q}=\frac{1}{qk} \left\{ \Sum{v=0}{r}\Sum{u=1}{q-1}-\frac{3}{2}(q-2u)^2x_{k-3u-v,q-u}x_{3u+v,u}\right.
\end{equation*}
\begin{equation*}
+3\Sum{v=0}{r}\Sum{u=0}{q-1}(m-3u-v)(m-3u-v-1)(m-3(q-u-1)+v-r)\cdot
\end{equation*}
\begin{equation*}
\cdot(m-3(q-u-1)+v-r-1)x_{3u+v,u}x_{k-3u-v-3,q-u-1} 
\end{equation*}
\begin{equation*}
  +\paren{(m-3q)qc_{r}x_{3q,q}-q(m-1)c_{1}x_{k-1,q}}\chi_{\llave{r\geq 1}}(r)
   \end{equation*}
\begin{equation*}
     \left.   +( (m-3q-1)qc_{1}x_{3q+1,q}-q(m-2)c_{2}x_{3q,q})\delta_{r,2}) \right\} . 
\end{equation*}
for every $k \in \natu$.
\end{teor}
\begin{demos}
Assume that for $0\leq l\leq k-1.$ If we write $k=3q+r$, $r\in \llave{0,1,2}$ and take $0\leq l \leq k-1$, $l=3u+v$, $v\in \llave{0,1,2}$ then 
$\deg_{t}(Q_{l})= u$ and note that $k-l=3q+r-(3u+v)=3(q-u)+r-v=3a+b$ with $b\in \llave{0,1,2}$ for $a=q-u$, $b=r-v$ if $r\geq v$ and $a=q-u-1$,
$b=3+r-v$ if $r<v$. Therefore, 
 \begin{equation*}
     \deg_{t}(Q_{k-l}) = \left\{
	       \begin{array}{ll}
		 q-u      & \mathrm{if\ } r\geq v \\
		 q-u-1 & \mathrm{if\ } r<v 
	       \end{array}
	     \right.
   \end{equation*}
and 
 \begin{equation*}
     \deg_{t}(Q_{k-l-3}) = \left\{
	       \begin{array}{ll}
		 q-u-1      & \mathrm{if\ } r\geq v \\
		 q-u-2 & \mathrm{if\ } r<v. 
	       \end{array}
	     \right.
   \end{equation*}
   
If we write $Q_{k}(m,t,c_{1},...,c_{k})=\Sum{j=0}{q}x_{k,j}(m,c_{1},...,c_{k})t^{j}$ then 
\begin{equation*}
x_{k,q}(m,c_{1},...,c_{k})t^{q}=
\frac{1}{k}\left\{ \Sum{l=1,l\equiv v (mod 3), v\geq r}{k-1}(m-l)\int W_{t}\paren{x_{k-l,q-u}(m,c_{1},...,c_{k-l})t^{q-u},
x_{l,u}(m,c_{1},...,c_{l})t^{u}}dt\right.
\end{equation*}
\begin{equation*}
  +3\Sum{l=0,l\equiv v (mod 3), v\geq r}{k-3}(m-l)(m-l-1)(m+l-k+3)(m+l-k+2)\cdot
   \end{equation*}
\begin{equation*}
  \left.
   \int \paren{x_{l,u}(m,c_{1},...,c_{l})t^{u}}
   \paren{x_{k-l-3,q-u-1}(m,c_{1},...,c_{k-l-3})t^{q-u-1}}dt
   \right\}.
\end{equation*}
Using the natural parametrization $l=3u+v$ and the simplified notation $x_{l,u}:=x_{l,u}(m,c_{1},...,c_{l})$ we obtain:
\begin{equation*}
x_{k,q}=\frac{1}{qk}\left\{ \Sum{v=0}{r}\Sum{\max\llave{0,1-v}\leq u \leq \min\llave{q,q+\frac{r-v-1}{3}}}{} (m-3u-v)(2u-q)x_{k-l,q-u}x_{l,u}\right.
\end{equation*}
\begin{equation*}
  \left. +3\Sum{v=0}{r}\Sum{u=0}{q-1}(m-3u-v)(m-3u-v-1)(m+3u+v-3q-r+3)(m+3u+v-3q-r+2)x_{l,u}x_{k-l-3,q-u-1}\right\}
\end{equation*}
\begin{equation*}
   = \frac{1}{qk}\left\{\Sum{1\leq u \leq q+\frac{r-1}{3}}{}(m-3u)(2u-q)x_{k-3u,q-u}x_{3u,u}\right.
\end{equation*}
\begin{equation*}
 +\Sum{1\leq u \leq q+\frac{r-2}{3}}{}(m-3u-1)(2u-q)x_{k-3u-1,q-u}x_{3u+1,u}\chi_{r\geq 1}(r)
\end{equation*}
\begin{equation*}
    +\Sum{u=1}{q-1}(m-3u-2)(2u-q)x_{k-3u-2,q-u}x_{3u+2,u}\delta_{r,2}
\end{equation*}
\begin{equation*}
  \left.  +3\Sum{v=0}{r}\Sum{u=0}{q-1}(m-3u-v)(m-3u-v-1)(m+3u+v-3q-r+3)(m+3u+v-3q-r+2)x_{l,u}x_{k-l-3,q-u-1}\right\}
\end{equation*}
\begin{equation*}
   = \frac{1}{qk}\left\{\Sum{u=1}{q-1}(m-3u)(2u-q)x_{k-3u,q-u}x_{3u,u}\right.
\end{equation*}
\begin{equation*}
    +(m-3q)qx_{r,0}x_{3q,q}\chi_{r\geq 1}(r)-q(m-1)x_{k-1,q}x_{1,0}\chi_{r\geq 1}(r)
\end{equation*}
\begin{equation*}
 +\Sum{u=1}{q-1}(m-3u-1)(2u-q)x_{k-3u-1,q-u}x_{3u+1,u}\chi_{r\geq 1}(r)
\end{equation*}
\begin{equation*}
    +(m-3q-1)qx_{1,0}x_{3q+1,q}\delta_{r,2}-q(m-2)x_{3q,q}x_{2,0}\delta_{r,2}
\end{equation*}
\begin{equation*}
    +\Sum{u=1}{q-1}(m-3u-2)(2u-q)x_{k-3u-2,q-u}x_{3u+2,u}\delta_{r,2}
\end{equation*}
\begin{equation*}
  \left.  +3\Sum{v=0}{r}\Sum{u=0}{q-1}(m-3u-v)(m-3u-v-1)(m+3u+v-3q-r+3)(m+3u+v-3q-r+2)x_{l,u}x_{k-l-3,q-u-1}\right\}
\end{equation*}
\begin{equation*}
 = \frac{1}{qk}\left\{ \Sum{v=0}{r}\Sum{u=1}{q-1} -\frac{3}{2}(q-2u)^2x_{k-3u-v,q-u}x_{3u+v,u}\right.   
\end{equation*}
\begin{equation*}
+3\Sum{v=0}{r}\Sum{u=0}{q-1}(m-3u-v)(m-3u-v-1)(m-3(q-u-1)+v-r)(m-3(q-u-1)+v-r-1)x_{l,u}x_{k-l-3,q-u-1}    
\end{equation*}
\begin{equation*}
    +(m-3q)qc_{r}x_{3q,q}\chi_{r\geq 1}(r)-q(m-1)c_{1}x_{k-1,q}\chi_{r\geq 1}(r)
\end{equation*}
\begin{equation*}
   \left. +\paren{(m-3q-1)qc_{1}x_{3q+1,q}-q(m-2)c_{2}x_{3q,q}}\delta_{r,2} \right\}.
\end{equation*}
Now, we write $x_{k,q}(m)=\Sum{j=0}{4q}a_{k,q,j}t^{j}$ and prove that the nice form of the  leading coefficient $a_{k,q,4q}$.

\begin{itemize}
    \item If $k\equiv 0$ (mod 3) then $r=0$ and the formula reduces to:
    \begin{equation*}
 x_{3q,q}= \frac{1}{3q^2}\left\{ \Sum{u=1}{q-1} -\frac{3}{2}(q-2u)^2x_{3(q-u),q-u}x_{3u,u}\right.   
\end{equation*}
\begin{equation*}
\left. +3\Sum{u=0}{q-1}(m-3u)(m-3u-1)(m-3(q-u-1))(m-3(q-u-1)-1)x_{3u,u}x_{3(q-u-1),q-u-1}    \right\}
\end{equation*}
   \begin{equation*}
=\frac{1}{q^2}\left\{ \Sum{u=1}{q-1} -\frac{1}{2}(q-2u)^2x_{3(q-u),q-u}x_{3u,u}\right.   
\end{equation*}
\begin{equation}\label{recur0}
\left. +3\Sum{u=0}{q-1}(m-3u)(m-3u-1)(m-3(q-u-1))(m-3(q-u-1)-1)x_{3u,u}x_{3(q-u-1),q-u-1}   \right\}. 
\end{equation}
Since $x_{0,0}(m)=1$ we have that $x_{0,0}\in \racio[m]$. If we assume that $x_{l,u}\in \racio[m]$, $l\equiv 0$ (mod 3), $0\leq l\leq k-1$ then \eqref{recur0} implies that $x_{k,q}\in \racio[m]$. Thus, \ref{1} follows by induction over $k$. 
We claim that $a_{3q,q,4q}=\frac{1}{q!}$. 

In fact, for $q=0$ we are okay, since $x_{0,0}=1$. Assume that $a_{3u,u,4u}=\frac{1}{u!}$ for $0\leq u \leq q-1$ and use \eqref{recur0} to obtain 
\begin{equation*}
    a_{3q,q,4q}=\frac{1}{q^2}\llave{-\Sum{u=1}{q-1}(q-2u)^2\frac{1}{(q-u)!}\frac{1}{u!}+\Sum{u=0}{q-1}\frac{1}{u!}\frac{1}{(q-u-1)!}}
    \end{equation*}
    \begin{equation*}
        =\frac{1}{q^2}\llave{\Sum{u=1}{q-1}u(q-u)\frac{1}{u!}\frac{1}{(q-u)!}-\Sum{u=1}{q-1}(q-u)^2\frac{1}{u!}\frac{1}{(q-u)!}
        +\Sum{u=0}{q-1}\frac{1}{u!}\frac{1}{(q-u-1)!}}
    \end{equation*}
    \begin{equation*}
        =\frac{1}{q^2}\llave{\frac{1}{(q-1)!}+\Sum{u=1}{q-1}(u+1)\frac{1}{u!}\frac{1}{(q-u-1)!}-\Sum{u=1}{q-1}(q-u)\frac{1}{u!}\frac{1}{(q-u-1)!}}
    \end{equation*}
    \begin{equation*}
      =\frac{1}{q^2}\llave{\frac{1}{(q-1)!}+\Sum{l=2}{q}l\frac{1}{(l-1)!}\frac{1}{(q-l)!}-\Sum{l=1}{q-1}l\frac{1}{(q-l)!}\frac{1}{(l-1)!}}  
    \end{equation*}
    \begin{equation*}
        =\frac{1}{q^2}\llave{\frac{1}{(q-1)!}+\frac{q}{(q-1)!}-\frac{1}{(q-1)!}}=\frac{1}{q!}.
    \end{equation*}
    
    \item If $k\equiv 1$ (mod 3) then $r=1$ and the formula reduces to:
    \begin{equation*}
        x_{3q+1,q}=\frac{1}{q(3q+1)}\left\{\Sum{w=0}{1}\Sum{u=1}{q-1}-\frac{3(q-2u)^2}{2}x_{3(q-u)+1-w,q-u}x_{3u+w,u}\right.
    \end{equation*}
    \begin{equation*}
       \left. +3\Sum{v=0}{1}\Sum{u=0}{q-1}(m-3u-v)(m-3u-v-1)\cdot \right.
           \end{equation*}
    \begin{equation*}
      \left.  \cdot(m-3(q-u-1)+v-1)(m-3(q-u-1)+v-2)x_{3u+v,u}x_{3(q-u-1)+(1-v),q-u-1} \right.
    \end{equation*}
    \begin{equation*}
     \left.   +(m-3q)qc_{1}x_{3q,q}-q(m-1)c_{1}x_{3q,q}\right\}
    \end{equation*}
    \begin{equation*}
        =\frac{1}{q(3q+1)}\left\{\Sum{u=1}{q-1}-3(q-2u)^2x_{3(q-u)+1,q-u}x_{3u,u}\right.
    \end{equation*}
       \begin{equation*}
        +6\Sum{v=0}{1}\Sum{u=0}{q-1}(m-3u)(m-3u-1)\cdot
           \end{equation*}
               \begin{equation}\label{recur1}
       \left. \cdot(m-3(q-u-1)-1)(m-3(q-u-1)-2)x_{3u,u}x_{3(q-u-1)+1,q-u-1}-q(3q-1)c_{1}x_{3q,q} \right\}.
    \end{equation}
    Since $x_{1,0}(m)=c_{1}$ we have that $x_{1,0}\in \racio[m]c_{1}$. If we assume that $x_{l,u}\in \racio[m]c_{1}$, $l\equiv 1$ (mod 3). Since 
    \eqref{recur1} tells us that $x_{k,q}$ is a sum of products  of the form $x_{l,u}x_{n,w}$ with $l\equiv 0$ (mod 3), $n\equiv 1$ (mod 3), then 
    $x_{l,u}\in \racio[m]$ and $x_{n,w}\in \racio[m]c_{1}$, in particular $x_{k,q}\in \racio[m]c_{1}.$ Thus, \ref{2} follows by induction over $k$. We claim that $a_{3q+1,q,4q}=\frac{c_{1}}{q!}=c_{1}a_{3q,q,4q}$.
    
    In fact, for $q=0$ we are okay, since $x_{1,0}=c_{1}$. Assume that $a_{3u+1,u,4u}=\frac{c_{1}}{u!}=c_{1}a_{3u,u,4u}$ for $0\leq u \leq q-1$ and use 
    \eqref{recur1} to obtain:
    \begin{equation*}
        a_{3q+1,q,4q}=\frac{1}{q(3q+1)}\left\{ \Sum{u=1}{q-1}-3(q-2u)^2 a_{3(q-u)+1,u,4u}a_{3u,u,4u}\right.
            \end{equation*}
       \begin{equation*}
       \left. +6\Sum{u=0}{q-1}a_{3u,u,4u}a_{3(q-u-1)+1,q-u-1,4(q-u-1)}-q(3q-1)c_{1}a_{3q,q,4q}\right\}
    \end{equation*}
    \begin{equation*}
        =\frac{1}{q(3q+1)}\left\{6c_{1}\left(\Sum{u=1}{q-1}-\frac{1}{2}(q-2u)^2 a_{3(q-u)+1,u,4u}a_{3u,u,4u}\right. \right.
          \end{equation*}
    \begin{equation*}
       \left. \left. +\Sum{u=0}{q-1}a_{3u,u,4u}a_{3(q-u-1)+1,q-u-1,4(q-u-1)}\right) -q(3q-1)c_{1}a_{3q,q,4q} \right\}
    \end{equation*}
    \begin{equation*}
        =\frac{1}{q(3q+1)}\llave{6c_{1}q^2a_{3q,q,4q}-q(3q-1)c_{1}a_{3q,q,4q}}=c_{1}a_{3q,q,4q}=\frac{c_{1}}{q!}.
    \end{equation*}

   \item If $k\equiv 2$ (mod 3) then $r=2$ and the formula reduces to:
   \begin{equation*}
       x_{3q+2,q}=\frac{1}{q(3q+2)}\left\{ \Sum{w=0}{2}\Sum{u=1}{q-1}-\frac{3}{2}(q-2u)^2x_{3(q-u)+2-w,q-u}x_{3u+w,u}\right.
   \end{equation*}
   \begin{equation*}
        \left. +3\Sum{v=0}{2}\Sum{u=0}{q-1}(m-3u-v)(m-3u-v-1)\cdot \right.
   \end{equation*}
  \begin{equation*}
      \left.  \cdot(m-3(q-u-1)+v-2)(m-3(q-u-1)+v-3)x_{3u+v,u}x_{3(q-u-1)+(2-v),q-u-1} \right.
    \end{equation*}
    \begin{equation*}
       \left. +(m-3q)qc_{2}x_{3q,q}-q(m-1)c_{1}x_{3q+1,q}+(m-3q-1)qc_{1}x_{3q+1,q}-q(m-2)c_{2}x_{3q,q}\right\}
    \end{equation*}
    \begin{equation*}
     =\frac{1}{q(3q+2)}\left\{ \Sum{u=1}{q-1}3(q-2u)^2x_{3(q-u)+2,q-u}x_{3u,u}\right.
   \end{equation*}
   \begin{equation*}
       +\Sum{u=1}{q-1}-\frac{3}{2}(q-2u)^2x_{3(q-u)+1,q-u}x_{3u+1,u}
   \end{equation*}
    \begin{equation*}
        \left. +6\Sum{u=0}{q-1}(m-3u)(m-3u-1)\cdot \right.
   \end{equation*}
  \begin{equation*}
      \left.  \cdot(m-3(q-u-1)-2)(m-3(q-u-1)-3)x_{3u,u}x_{3(q-u-1)+2,q-u-1} \right.
    \end{equation*}  
   \begin{equation*}
        \left. +3\Sum{u=0}{q-1}(m-3u-1)(m-3u-2)\cdot \right.
   \end{equation*}
  \begin{equation*}
      \left.  \cdot(m-3(q-u-1)-1)(m-3(q-u-1)-1)x_{3u+1,u}x_{3(q-u-1)+1,q-u-1} \right.
    \end{equation*}
    \begin{equation}\label{recur2}
       \left. -(3q-2)qc_{2}x_{3q,q}-3q^2c_{1}x_{3q+1,q}\right\}.
    \end{equation}
 Since $x_{2,0}(m)=c_{2}$ we have that $x_{2,0}\in \racio[m]c_{1}^2+\racio[m]c_{2}$. If we assume that $x_{l,u}\in \racio[m]c_{1}^2+\racio[m]c_{2}$,
 $l\equiv 2$ (mod 3). Since \eqref{recur2} tells us that $x_{k,q}$ is a sum of products of the form $x_{l,u}x_{n,w}$ with $l\equiv 0$ (mod 3), 
 $n\equiv 2$ (mod 3) or $l\equiv 1$ (mod 3), $n\equiv 1$ (mod 3) we have that $x_{l,u}\in \racio[m]$ and $x_{n,w}\in \racio[m]c_{1}^2+\racio[m]c_{2}$ 
 or $x_{l,u},x_{n,w}\in \racio[m]c_{1}$. Therefore, $x_{l,u},x_{n,w}\in \racio[m]c_{1}^2+\racio[m]c_{2}.$ Thus, \ref{3} follows by induction over $k$.
 
 We claim that $a_{3q+2,q,4q}=\frac{c_{2}}{q!}=c_{2}a_{3q,q,4q}$. In fact, if $q=0$ we are okay. Since, $x_{2,0}=c_{2}.$ 
 Assume that $a_{3u+2,u,4u}=\frac{c_{2}}{u!}=c_{2}a_{3u,u,4u}$, for $0\leq u\leq q-1$ and we use \eqref{recur2} to obtain 
 \begin{equation*}
     a_{3q+2,q,4q}=\frac{1}{q(3q+2)}\left\{\Sum{u=1}{q-1}-3(q-2u)^2a_{3(q-u)+2,q-u,4(q-u)}a_{3u,u,4u}\right.
    \end{equation*}
    \begin{equation*}
       +\left.\Sum{u=1}{q-1}-\frac{3}{2}(q-2u)^2 a_{3(q-u)+1,q-u,4(q-u)}a_{3u+1,u,4u}\right.
    \end{equation*}
    \begin{equation*}
        +6\Sum{u=0}{q-1}a_{3u,u,4u}a_{3(q-u-1)+2,q-u-1,4(q-u-1)}
    \end{equation*}
    \begin{equation*}
       \left. +3\Sum{u=0}{q-1}a_{3u+1,u,4u}a_{3(q-u-1)+1,q-u-1,4(q-u-1)}\right.
    \end{equation*}
    \begin{equation*}
       \left. -(3q-2)qc_{2}a_{3q,q,4q}-3q^2c_{1}a_{3q+1,q,4q}\right\}
    \end{equation*}
    \begin{equation*}
        =\frac{1}{q(3q+2)}\left\{6c_{2}\left(\Sum{u=1}{q-1}-\frac{1}{2}(q-2u)^2 a_{3(q-u),u,4u} a_{3u,u,4u}\right.\right.
    \end{equation*}
    \begin{equation*}
        \left.\left. +\Sum{u=0}{q-1}a_{3u,u,4u}a_{3(q-u-1),q-u-1,4(q-u-1)} \right)\right.
    \end{equation*}
    \begin{equation*}
        +3c_{1}\left(\Sum{u=1}{q-1}-\frac{1}{2}(q-2u)^2 a_{3(q-u),u,4u} a_{3u,u,4u}\right.
    \end{equation*}
    \begin{equation*}
       \left.\left. +\Sum{u=0}{q-1}a_{3u,u,4u}a_{3(q-u-1),q-u-1,4(q-u-1)} \right)\right. 
    \end{equation*}
    \begin{equation*}
       \left. -(3q-2)qc_{2}a_{3q,q,4q}-3q^2c_{1}^2 a_{3q,q,4q}\right\}
    \end{equation*}
    \begin{equation*}
        =\frac{1}{q(3q+2)}\left\{6c_{2}q^2 a_{3q,q,4q}+3c_{1}q^2 a_{3q,q,4q}-(3q-2)qc_{2}a_{3q,q,4q}-3q^2 c_{1}^2a_{3q,q,4q}\right\}
    \end{equation*}
    \begin{equation*}
        =c_{2}a_{3q,q,4q}=\frac{c_{2}}{q!}.
    \end{equation*} \QED
    \end{itemize}
\end{demos}
Note that the only formula for the leading polynomial $x_{k,q}$ which depends on elements of the same class is in the case $k\equiv 0$ (mod 3) and is the simplest. We will use that to show that $Rf\neq 0$ for $m\geq 5$.

\begin{coro}
We have the identities 
\begin{equation*}
    \deg_{t}(Q_{k})=q \hspace{0.1cm} \text{and} \hspace{0.1cm} \deg_{m}(Q_{k})=4q. 
\end{equation*}
\end{coro}
\begin{demos}
We have the recursion $Q_{0}=1$, $Q_{1}=c_{1}$, $Q_{2}=c_{2}$ and by Proposition \ref{recurQ} 
    \begin{equation*}
    Q_{k}(t)=\frac{1}{k}\left\{\Sum{l=1}{k-1}(m-l)\int W_{t}(Q_{k-l}(t),Q_{l}(t))dt\right.
\end{equation*}   
\begin{equation*}
    \left.+3\Sum{l=0}{k-3}(m-l)(m-l-1)(m+l-k+3)(m+l-k+2)\int Q_{l}(t)Q_{k-l-3}(t)dt\right\}+c_{k}
\end{equation*}
for an arbitrary $c_{k}\in \complex.$
We take the symmetrization in the first sum 
\begin{equation*}
    \Sum{l=1}{k-1}(m-l)\int W_{t}(Q_{k-l}(t),Q_{l}(t))dt=
    \Sum{l=1}{\corch{\frac{k-1}{2}}}(m-l)\int W_{t}(Q_{k-l}(t),Q_{l}(t))dt
    \end{equation*}
\begin{equation*}
+\Sum{l=\corch{\frac{k-1}{2}}+1}{k-1}(m-l)\int W_{t}(Q_{k-l}(t),Q_{l}(t))dt
 =\Sum{l=1}{\corch{\frac{k-1}{2}}}(m-l)\int W_{t}(Q_{k-l}(t),Q_{l}(t))dt
 \end{equation*}
\begin{equation*}
    - \Sum{l=1}{\corch{\frac{k-1}{2}}}(m-(k-l))\int W_{t}(Q_{k-l}(t),Q_{l}(t))dt
    =\Sum{l=1}{k-1}\frac{1}{2}(k-2l)\int W_{t}(Q_{k-l}(t),Q_{l}(t))dt.
\end{equation*}
To obtain 
   \begin{equation*}
    Q_{k}(t)=\frac{1}{k}\left\{\Sum{l=1}{k-1}\frac{1}{2}(k-2l)\int W_{t}(Q_{k-l}(t),Q_{l}(t))dt\right.
\end{equation*}   
\begin{equation*}
    \left.+3\Sum{l=0}{k-3}(m-l)(m-l-1)(m+l-k+3)(m+l-k+2)\int Q_{l}(t)Q_{k-l-3}(t)dt\right\}+c_{k}
\end{equation*}
for an arbitrary $c_{k}\in \complex.$
Therefore, we find by induction over $q$ that $\deg_{m}(Q_{k})\leq 4q$.
On the other hand, the previous theorem implies $x_{k,q}$ has leading coefficient $a_{k,q,4q}=\frac{c_{r}}{q!}\neq 0$, therefore $\deg_{m} (Q_{k})\geq \deg_{m} (x_{k,q})=4q.$ Furthermore, $\deg_{t}(Q_{k})= q$. \QED
\end{demos}

\section{The Operator Remainder as a Function of the Polynomials}\label{sec3}

In this section we evaluate the operator remainder $R$ in the fundamental polynomial $f$ such that $Bf=0$ to obtain a family of polynomials in the $Q_{k}$.

Remember the Equation \eqref{remainder} and note that
\begin{equation*}
    Rf(x,t)=\Sum{p=0}{m-2}\llave{\Sum{j=0}{m}jW(f_{j},f_{p+1-j})(t)
    +3\Sum{j=0}{m}j(j-1)(p+4-j)(p+3-j)f_{j}(t)f_{p+4-j}(t)}x^p
\end{equation*}
\begin{equation*}
    =f_{m}(t)^2\Sum{p=0}{m-2}\left\{\Sum{l=0}{m}(m-l)W(P_{m,m-l},P_{m,p+1-m+l})(t)\right.
    \end{equation*}
\begin{equation*}
    \left.+3\Sum{l=0}{m}(m-l)(m-l-1)(p+4-m+l)(p+3-m+l)P_{m,m-l}(t)P_{m,p+4-m+l}(t)\right\}x^p
\end{equation*}
\begin{equation*}
    =f_{m}(t)^2\Sum{k=m+1}{2m-1}\left\{\Sum{l=0}{m}(m-l)W(P_{m,m-l},P_{m,m-(k-l)})(t)\right.
    \end{equation*}
\begin{equation*}
      \left.+3\Sum{l=0}{m}(m-l)(m-l-1)(m+l-k+3)(m+l-k+2)P_{m,m-l}(t)P_{m,m-(k-l-3)}(t)\right\}x^{2m-k-1}
\end{equation*}
\begin{equation*}
    =f_{m}(t)^2\Sum{k=m+1}{2m-1}\left\{\Sum{l=0}{m}(m-l)W_{t}(Q_{l}(t),Q_{k-l}(t))\right.
    \end{equation*}
\begin{equation*}
      \left.+3\Sum{l=0}{m}(m-l)(m-l-1)(m+l-k+3)(m+l-k+2)Q_{l}(t)Q_{k-l-3}(t)\right\}x^{2m-k-1}
\end{equation*}
\begin{equation*}
    =f_{m}(t)^2\Sum{k=m+1}{2m-1}R_{k}(m,t)x^{2m-k-1}
\end{equation*}
with 
\begin{equation*}
    R_{k}(t)=\Sum{l=k-m}{m}(m-l)W_{t}(Q_{l}(t),Q_{k-l}(t))
\end{equation*}
\begin{equation*}
    +3 \Sum{l=k-m-1}{m-2}(m-l)(m-l-1)(m+l-k+3)(m+l-k+2)Q_{l}(t)Q_{k-l-3}(t),
\end{equation*}
for every $m+1\leq k \leq 2m-1$.

However, we can symmetrize the first term 
\begin{equation*}
    \Sum{l=k-m}{m}(m-l)W_{t}(Q_{l}(t),Q_{k-l}(t))
    =\Sum{l=k-m}{\corch{\frac{k-1}{2}}}(m-l)W_{t}(Q_{l}(t),Q_{k-l}(t))
    +\Sum{l=\corch{\frac{k-1}{2}}+1}{m}(m-l)W_{t}(Q_{l}(t),Q_{k-l}(t))
\end{equation*}
\begin{equation*}
    =\Sum{l=k-m}{\corch{\frac{k-1}{2}}}(m-l)W_{t}(Q_{l}(t),Q_{k-l}(t))
    -\Sum{j=k-m}{\corch{\frac{k-1}{2}}}(m-(k-j))W_{t}(Q_{l}(t),Q_{k-j}(t))
\end{equation*}
\begin{equation*}
    =\Sum{l=k-m}{\corch{\frac{k-1}{2}}}(k-2l)W_{t}(Q_{l}(t),Q_{k-l}(t))
    =\Sum{l=k-m}{m}\frac{(k-2l)}{2} W_{t}(Q_{l}(t),Q_{k-l}(t))
\end{equation*}
Therefore, 
\begin{equation*}
    R_{k}(t)=\Sum{l=k-m}{m}\frac{(k-2l)}{2} W_{t}(Q_{l}(t),Q_{k-l}(t))
\end{equation*}
\begin{equation*}
    +3 \Sum{l=k-m-1}{m-2}(m-l)(m-l-1)(m+l-k+3)(m+l-k+2)Q_{l}(t)Q_{k-l-3}(t).
\end{equation*}

\subsection{The Leading Function of the Remainder Polynomials}\label{sub3.1}

Let $k=3q+r$ for $r\in \llave{0,1,2}$ and write $R_{k}(m,t)=\Sum{j=0}{q-1}y_{k,j}(m)t^j$. Using the parametrization $l=3u+v$ for $v\in \llave{0,1,2}$ we have 
\begin{equation*}
    y_{k,q-1}(m)t^{q}=\Sum{l=k-m,l\equiv v (mod 3), r\geq v}{m} \frac{(k-2l)}{2} W_{t}(x_{l,u}t^u , x_{k-l,q-u}t^{q-u})
\end{equation*}
\begin{equation*}
    +3\Sum{l=k-m-1,l\equiv v (mod 3), r\geq v}{m-2}(m-l)(m-l-1)(m+l-k+3)(m+l-k+2)(x_{l,u}t^u)(x_{k-l-3,q-u-1}t^{q-u-1})
\end{equation*}
\begin{equation*}
    =\left\{\Sum{l=k-m,l\equiv v (mod 3), r\geq v}{m} \frac{(k-2l)}{2} (q-2u)x_{l,u}x_{k-l,q-u}\right.
\end{equation*}
\begin{equation*}
    \left.+3\Sum{l=k-m-1,l\equiv v (mod 3), r\geq v}{m-2}(m-l)(m-l-1)(m+l-k+3)(m+l-k+2)x_{l,u}x_{k-l-3,q-u-1}\right\}t^{q-1}.
\end{equation*}
Therefore, 
\begin{equation*}
    y_{k,q-1}(m)=\Sum{l=k-m,l\equiv v (mod 3), r\geq v}{m} \frac{(k-2l)}{2} (q-2u)x_{l,u}x_{k-l,q-u}
\end{equation*}
\begin{equation*}
    +3\Sum{l=k-m-1,l\equiv v (mod 3), r\geq v}{m-2}(m-l)(m-l-1)(m+l-k+3)(m+l-k+2)x_{l,u}x_{k-l-3,q-u-1}.
\end{equation*}

In particular, if $k\equiv 0$ (mod 3) then $r=0$ and 
    \begin{equation*}
    y_{k,q-1}(m)=\Sum{l=k-m,l\equiv 0 (mod 3)}{m} \frac{3}{2}(q-2u)^2 x_{l,u}x_{k-l,q-u}
\end{equation*}
\begin{equation*}
    +3\Sum{l=k-m-1,l\equiv 0 (mod 3)}{m-2}(m-l)(m-l-1)(m+l-k+3)(m+l-k+2)x_{l,u}x_{k-l-3,q-u-1}.
\end{equation*}

 From now on we denote $z_{q}=x_{3q,q}$ and $Y_{q}=y_{3q,q-1}$, for $q\in \natu$.
 \begin{remark}
  Note that in a first glance the functions $y_{k,q-1}(m)$, $m+1\leq k \leq 2m-1$ are not polynomials in the variable $m$, since the sum defining it depends on $m$, although they are for $k$ near $2m-1$. However, we shall use that they are elements of $\racio[x_{l,u}\mid 0\leq l \leq m]$ to understand its simultaneous zeros in terms of the algebraic sets defined by the polynomials $x_{l,u}, 0\leq l \leq m$.
 \end{remark}
 
 \subsection{Algebraic Sets defined by the Leading Functions of the Remainder Polynomials}\label{sub3.2}

The following theorem gives us a characterization of the zero set defined by the functions $Y_{q}$, $m+1\leq k\leq 2m-1$ with $k=3q+r$, $r\in \llave{0,1,2}$ in terms of the algebraic sets of $z_{q}=x_{3q,q}$.

\begin{lemma}
\begin{enumerate}
\item If $m\equiv 0$ (mod 3) then
\begin{equation*}
    Y_{q}=\Sum{u=q-s}{s}\frac{3}{2}(q-2u)^2 z_{u}z_{q-u}
    \end{equation*}
\begin{equation*}
    +27 \Sum{u=q-s}{s-1}(s-u)(3s-3u-1)(s+u-q+1)(3s+3u-3q+2)z_{u}z_{q-u-1},
\end{equation*}
with $m=3s$ and $s+1\leq q\leq 2s-1$.

Furthermore, for $s\geq 2$ we have
\begin{equation*}
    V(Y_{2s-1})=V(z_{s}+36z_{s-1})\cup V(z_{s-1})
\end{equation*}
and 
\begin{equation*}
    \bigcap_{q=s+v}^{2s-1}V(Y_{q})=\corch{V(z_{s}+36z_{s-1}) \cap \bigcap_{u=v}^{s-2}V(z_{u})}\cup \corch{\bigcap_{u=v}^{s-1}V(z_{u})}
    \cup \corch{\bigcap_{u=v+1}^{s}V(z_{u})},
\end{equation*}
for every $1\leq v \leq s-2$. 

\item If $m\equiv 1$ (mod 3) then
\begin{equation*}
    Y_{q}=\Sum{u=q-s}{s}\frac{3}{2}(q-2u)^2  z_{u}z_{q-u}
\end{equation*}
\begin{equation*}
    +27 \Sum{u=q-s}{s-1}(s-u)(3(s-u)+1)(s+u-q+1)(3s+3u-3q+4)z_{u}z_{q-u-1},
\end{equation*}
with $m=3s+1$ and $s+1\leq q \leq 2s$.

Furthermore, for $s\geq 2$ we have
\begin{equation*}
    V(Y_{2s-1})=V(z_{s}+144z_{s-1})\cup V(z_{s-1})
\end{equation*}
and 
\begin{equation*}
    \bigcap_{q=s+v}^{2s-1}V(Y_{q})=\corch{V(z_{s}+144z_{s-1}) \cap \bigcap_{u=v}^{s-2}V(z_{u})}\cup \corch{\bigcap_{u=v}^{s-1}V(z_{u})}
    \cup \corch{\bigcap_{u=v+1}^{s}V(z_{u})},
\end{equation*}
for every $1\leq v \leq s-2$. 

    \item If $m\equiv 2$ (mod 3) then
    \begin{equation*}
        Y_{q}=\Sum{u=q-s}{s}\frac{3}{2}(q-2u)^2  z_{u}z_{q-u}
    \end{equation*}
    \begin{equation*}
        +3\Sum{u=q-s-1}{s}(3(s-u)+2)(3(s-u)+1)(3(s-q+u+1)+2)(3(s-q+u+1)+1)z_{u}z_{q-u-1},
    \end{equation*}
with $m=3s+2$ and $s+1\leq q \leq 2s+1$.

Furthermore, for $s\geq 1$ we have

\begin{equation*}
    \bigcap_{l=v}^{s}V(Y_{2l+1})=\bigcap_{u=v}^{s}V(z_{u}),\hspace{0.1cm} \text{for every} \hspace{0.1cm} 1\leq v\leq s.
\end{equation*}
\end{enumerate}
\end{lemma}
\begin{demos}
\begin{enumerate}
 \item Suppose that $m\equiv 0$ (mod 3).
 
 Since $Y_{2s-1}=3z_{s-1}(z_{s}+36z_{s-1})$ we have that $V(Y_{2s-1})=V(z_{s}+36z_{s-1})\cup V(z_{s-1})$ and we obtain the result for $v=s$.
 On the other hand, $Y_{2s-2}=12z_{s-2}(z_{s}+90z_{s-1})$ implies that:
 \begin{equation*}
     \bigcap_{q=2s-2}^{2s-1}V(Y_{q})=\corch{V(z_{s}+36z_{s-1})\cup V(z_{s-1})}\cap 
     \corch{V(z_{s-2})\cup V(z_{s}+90z_{s-1})}
 \end{equation*}
 \begin{equation*}
     =\corch{V(z_{s}+36z_{s-1})\cap V(z_{s-2})}\cap \corch{\bigcap_{u=s-2}^{s-1}V(z_{u})}\cup \corch{\bigcap_{u=s-1}^{s}V(z_{u})}
 \end{equation*}
 Assume that 
 \begin{equation*}
    \bigcap_{q=s+v+1}^{2s-1}V(Y_{q})=\corch{V(z_{s}+36z_{s-1}) \cap \bigcap_{u=v+1}^{s-2}V(z_{u})}\cup \corch{\bigcap_{u=v+1}^{s-1}V(z_{u})}
    \cup \corch{\bigcap_{u=v+2}^{s}V(z_{u})}
\end{equation*}
 and note that 
 \begin{equation*}
     Y_{s+v}=\Sum{u=v}{s}\frac{3}{2}(s+v-2u)^2z_{u}z_{s+v-u}
 \end{equation*}
 \begin{equation*}
     +27\Sum{u=v}{s-1}(s-u)(3s-3u-1)(u-v+1)(3(u-v)+2)z_{u}z_{s+v-u-1}.
 \end{equation*}
 
 Then,
\begin{equation*}
    \bigcap_{q=s+v}^{2s-1}V(Y_{q})=V(Y_{s+v})\cap  \bigcap_{q=s+v+1}^{2s-1}V(Y_{q})
     \end{equation*}
 \begin{equation*}
  =V(Y_{s+v})\cap \left\{ \corch{V(z_{s}+36z_{s-1})\cap \bigcap_{u=v+1}^{s-2}V(z_{u})}
   \cup \corch{\bigcap_{u=v+1}^{s-1}V(z_{u})} \cup \bigcap_{u=v+2}^{s}V(z_{u}) \right\}
\end{equation*}
\begin{equation*}
    =\left[ V(3(s-v)^2z_{v}z_{s}+108(s-v)(3s-3v-1)z_{v}z_{s-1})\cap V(z_{s}+36z_{s-1}) \cap \bigcap_{u=v+1}^{s-2} V(z_{u}) \right]
\end{equation*}
\begin{equation*}
    \cup \corch{V(3(s-v)^2z_{v}z_{s})\cap \bigcap_{u=v+1}^{s-1}V(z_{u})} \cup \corch{\bigcap_{u=v+2}^{s}V(z_{u})}
\end{equation*}
\begin{equation*}
    =\left[ V((s-v)z_{v}((s-v)z_{s}+36(3s-3v-1)z_{s-1})\cap V(z_{s}+36z_{s-1}) \cap \bigcap_{u=v+1}^{s-2} V(z_{u}) \right]
\end{equation*}
\begin{equation*}
    \cup \corch{\bigcap_{u=v}^{s-1} V(z_{u})} \cup \corch{ \bigcap_{u=v+1}^{s} V(z_{u})}
\end{equation*}
\begin{equation*}
    =\corch{V(z_{s}+36z_{s-1})\cap \bigcap_{u=v}^{s-2} V(z_{u}) } \cup \corch{\bigcap_{u=v}^{s-1} V(z_{u})} \cup \corch{\bigcap_{u=v+1}^{s} V(z_{u})}
\end{equation*}
Thus, the assertion follows by induction.

    \item Suppose that $m\equiv 1$ (mod 3).
    
    Since $Y_{2s-1}=3z_{s-1}(z_{s}+144z_{s-1})$ we have that $V(Y_{2s-1})=V(z_{s}+144z_{s-1})\cup V(z_{s-1})$ and we obtain the result for $v=s$.
 On the other hand, $Y_{2s-2}=12z_{s-2}(z_{s}+252z_{s-1})$ implies that:
 \begin{equation*}
     \bigcap_{q=2s-2}^{2s-1}V(Y_{q})=\corch{V(z_{s}+144z_{s-1})\cup V(z_{s-1})}\cap 
     \corch{V(z_{s-2})\cup V(z_{s}+252z_{s-1})}
 \end{equation*}
 \begin{equation*}
     =\corch{V(z_{s}+144z_{s-1})\cap V(z_{s-2})}\cap \corch{\bigcap_{u=s-2}^{s-1}V(z_{u})}\cup \corch{\bigcap_{u=s-1}^{s}V(z_{u})}
 \end{equation*}
 Assume that 
 \begin{equation*}
    \bigcap_{q=s+v+1}^{2s-1}V(Y_{q})=\corch{V(z_{s}+144z_{s-1}) \cap \bigcap_{u=v+1}^{s-2}V(z_{u})}\cup \corch{\bigcap_{u=v+1}^{s-1}V(z_{u})}
    \cup \corch{\bigcap_{u=v+2}^{s}V(z_{u})}
\end{equation*}
 and note that 
 \begin{equation*}
     Y_{s+v}=\Sum{u=v}{s}\frac{3}{2}(s+v-2u)^2z_{u}z_{s+v-u}
 \end{equation*}
 \begin{equation*}
     +27\Sum{u=v}{s-1}(s-u)(3(s-u)-1)(u-v+1)(3(u-v+1)+1)z_{u}z_{s+v-u-1}.
 \end{equation*}
 
 Then,
\begin{equation*}
    \bigcap_{q=s+v}^{2s-1}V(Y_{q})=V(Y_{s+v})\cap  \bigcap_{q=s+v+1}^{2s-1}V(Y_{q})
     \end{equation*}
 \begin{equation*}
  =V(Y_{s+v})\cap \left\{ \corch{V(z_{s}+144z_{s-1})\cap \bigcap_{u=v+1}^{s-2}V(z_{u})}
   \cup \corch{\bigcap_{u=v+1}^{s-1}V(z_{u})} \cup \bigcap_{u=v+2}^{s}V(z_{u}) \right\}
\end{equation*}
\begin{equation*}
    =\left[ V(3(s-v)^2z_{v}z_{s}+216(s-v)(3s-3v-1)z_{v}z_{s-1})\cap V(z_{s}+144z_{s-1}) \cap \bigcap_{u=v+1}^{s-2} V(z_{u}) \right]
\end{equation*}
\begin{equation*}
    \cup \corch{V(3(s-v)^2z_{v}z_{s})\cap \bigcap_{u=v+1}^{s-1}V(z_{u})} \cup \corch{\bigcap_{u=v+2}^{s}V(z_{u})}
\end{equation*}
\begin{equation*}
    =\left[ V((s-v)z_{v}((s-v)z_{s}+72(3s-3v-1)z_{s-1})\cap V(z_{s}+144z_{s-1}) \cap \bigcap_{u=v+1}^{s-2} V(z_{u}) \right]
\end{equation*}
\begin{equation*}
    \cup \corch{\bigcap_{u=v}^{s-1} V(z_{u})} \cup \corch{ \bigcap_{u=v+1}^{s} V(z_{u})}
\end{equation*}
\begin{equation*}
    =\corch{V(z_{s}+144z_{s-1})\cap \bigcap_{u=v}^{s-2} V(z_{u}) } \cup \corch{\bigcap_{u=v}^{s-1} V(z_{u})} \cup \corch{\bigcap_{u=v+1}^{s} V(z_{u})}
\end{equation*}
Thus, the assertion follows by induction.

    \item Suppose that $m\equiv 2$ (mod 3).
    
    Since $Y_{2s+1}=z_{s}^2$ we have that $V(Y_{2s+1})=V(z_{s})$ and we obtain the result for $v=s$. Assume that $ \bigcap_{l=v+1}^{s}V(Y_{2l+1})=\bigcap_{u=v+1}^{s}V(z_{u})$ and note that 
\begin{equation*}
Y_{2v+1}=\Sum{u=2v+1-s}{s}  \frac{3}{2}(2(v-u)+1)^2  z_{u}z_{2v+1-u}
   \end{equation*} 
   \begin{equation*}
+3 \Sum{u=2v-s}{s}(3(s-u)+2)(3(s-u)+1)(3(s-2v+u)+2)(3(s-2v+u)+1)z_{u}z_{2v-u}.
\end{equation*}
Then,
\begin{equation*}
    \bigcap_{l=v}^{s}V(Y_{2l+1})=V(Y_{2v+1})\cap \bigcap_{l=v+1}^{s}V(Y_{2l+1})
    =V(Y_{2v+1})\cap \bigcap_{u=v+1}^{s}V(z_{u})
   \end{equation*} 
   \begin{equation*}
    =V(3(3(s-v)+2)^2(3(s-v)+1)^2z_{v}^2)\cap \bigcap_{u=v+1}^{s}V(z_{u})
    =\bigcap_{u=v}^{s}V(z_{u}).
    \end{equation*}
    Thus, the assertion follows by induction. \QED
\end{enumerate}
\end{demos}

\begin{remark}
 In particular, we have that the irreducible components of the intersection of hypersurfaces defined by $Y_{q}$, $s+1\leq q\leq 2s-1$ are hyperplanes defined by the coordinates $z_{u}$, $1\leq u \leq s$. 
\end{remark}

The following corollary is the algebraic counterpart of the previous theorem.

\begin{coro}
\begin{enumerate}
     \item If $m\equiv 0$ (mod 3), $m=3s$ then for every $s\geq 2$
     $$\langle Y_{2s-1} \rangle =\langle z_{s}+36z_{s-1}\rangle \cap \langle z_{s-1}\rangle $$ 
     and 
     $$\sqrt{\langle Y_{q}\mid s+v\leq q\leq 2s-1\rangle} = \langle z_{s}+36z_{s-1}, z_{u},v\leq u \leq s-2\rangle $$
   $$  \cap \langle z_{u}\mid v\leq u \leq s-1\rangle \cap \langle z_{u}\mid v+1\leq u \leq s \rangle,$$
     for every $1\leq v\leq s-2$.
     
    \item If $m\equiv 1$ (mod 3), $m=3s+1$ then for every $s\geq 2$
   $$\langle Y_{2s-1} \rangle =\langle z_{s}+144z_{s-1}\rangle \cap \langle z_{s-1}\rangle $$ 
     and 
     $$\sqrt{\langle Y_{q}\mid s+v\leq q\leq 2s-1\rangle} = \langle z_{s}+144z_{s-1}, z_{u},v\leq u \leq s-2\rangle $$
   $$  \cap \langle z_{u}\mid v\leq u \leq s-1\rangle \cap \langle z_{u}\mid v+1\leq u \leq s \rangle,$$
     for every $1\leq v\leq s-2$.
     
    \item If $m\equiv 2$ (mod 3), , $m=3s+2$ then for every $s\geq 1$
$$\sqrt{\langle Y_{2l+1}\mid v\leq l \leq s\rangle}= \langle z_{u}\mid v\leq u \leq s\rangle,$$
for every $1\leq v\leq s$.
\end{enumerate}

\end{coro}
\begin{demos}
The proof is a direct consequence of the theorem and the Hilbert Nullstellensatz.
\end{demos}

The following corollary will be used to characterize the solutions of Equation \eqref{bildiff}.

\begin{coro}\label{leadrestric}
\begin{enumerate}
 
\item If $m\equiv 0$ (mod 3) then for every $s\geq 2$
\begin{equation*}
    V(Y_{2s-1})=V(z_{s}+36z_{s-1})\cup V(z_{s-1})
\end{equation*}
and 
\begin{equation*}
    \bigcap_{q=s+1}^{2s-1}V(Y_{q})=\corch{V(z_{s}+36z_{s-1}) \cap \bigcap_{u=1}^{s-2}V(z_{u})}\cup \corch{\bigcap_{u=1}^{s-1}V(z_{u})}
    \cup \corch{\bigcap_{u=2}^{s}V(z_{u})},
\end{equation*}
with $m=3s$.

\item If $m\equiv 1$ (mod 3) then for every $s\geq 2$
\begin{equation*}
    V(Y_{2s-1})=V(z_{s}+144z_{s-1})\cup V(z_{s-1})
\end{equation*}
and 
\begin{equation*}
    \bigcap_{q=s+1}^{2s-1}V(Y_{q})=\corch{V(z_{s}+144z_{s-1}) \cap \bigcap_{u=1}^{s-2}V(z_{u})}\cup \corch{\bigcap_{u=1}^{s-1}V(z_{u})}
    \cup \corch{\bigcap_{u=2}^{s}V(z_{u})},
\end{equation*}
with $m=3s+1$.
    \item If $m\equiv 2$ (mod 3) then for every $s\geq 1$
\begin{equation*}
    \bigcap_{l=1}^{s}V(Y_{2l+1})=\bigcap_{u=1}^{s}V(z_{u}),
\end{equation*}
with $m=3s+2$.
\end{enumerate}
\end{coro}

\section{Characterization of Solutions of the Bilinear Differential and KdV-like equations}\label{sec4}

We conclude with the expected characterization theorem for polynomial solutions of Equation \eqref{bildiff} and their associated rational solutions of Equation \eqref{KdVlike}.

\begin{teor}
The only fundamental solutions of the  bilinear differential equation
\begin{equation*}
        2f_{xt}f-2f_{t}f_{x}+6(f_{xx})^2=0
\end{equation*}
are 
\begin{equation*}
f(x,t)=1,   
\end{equation*}
\begin{equation*}
f(x,t)=x+c_{0}, 
\end{equation*}
\begin{equation*}
f(x,t)=x^3+c_{2}x^2+\frac{1}{3}c_{2}^2x+36t+c_{0}    
\end{equation*}
and
\begin{equation*}
f(x,t)=x^4+c_{3}x^3+c_{2}x^2+(144t+c_{1})x+36c_{3}t+\frac{1}{4}c_{1}c_{3}-\frac{1}{12}c_{2}^2
\end{equation*}
for arbitrary constants $c_{0},c_{1},c_{2},c_{3}$.
\end{teor}
\begin{demos}
By simple calculations we have that $f(x,t)=1$, $f(x,t)=x+c_{0}$ are solutions. We don't have solutions of degree $m=2$, since if 
$f(x,t)=x^2+c_{1}x+c_{0}$ then $Tf=12\neq 0$.

If we consider the case $m=3$ then $P_{3,3}(t)=1$, $P_{3,2}(t)=c_{2}$, $P_{3,1}(t)=c_{1}$ and
\begin{equation*}
    P_{3,0}(t)= \frac{1}{3}\int \llave{\frac{3}{16}(6)^2(4)^2 P_{3,3}(t)^2+W(P_{3,2},P_{3,1})(t)}dt
    =36\int dt =36t+c_{0}.
\end{equation*}
The only restriction is 
\begin{equation*}
    W(P_{3,1},P_{3,0})(t)+3\Sum{j=1}{3}j(j-1)(4-j)(3-j)P_{3,j}(t)P_{3,4-j}(t)
    =-36c_{1}+12c_{2}^2=0.
\end{equation*}
Hence, $c_{1}=\frac{1}{3}c_{2}^2.$

Therefore, the fundamental solution is 
\begin{equation*}
f(x,t)=x^3+c_{2}x^2+\frac{1}{3}c_{2}^2x+36t+c_{0}.    
\end{equation*}

In the case $m=4$ we have $P_{4,4}(t)=1$, $P_{4,3}(t)=c_{3}$, $P_{4,2}(t)=c_{2}$ and 
\begin{equation*}
    P_{4,1}(t)=\frac{1}{3}\int\llave{\frac{3}{16}(8)^2(6)^2P_{4,4}(t)^2+W(P_{4,3},P_{4,2})(t)}dt
    =144 \int dt=144t+c_{1},
\end{equation*}

\begin{equation*}
    P_{4,0}(t)=\frac{1}{4}\int \llave{6(3)(4)(3)P_{4,3}(t)P_{4,4}(t)+2W(P_{4,3},P_{4,1})(t)}dt
  \end{equation*}  
  \begin{equation*}
    =108\int c_{3}dt-\frac{1}{2}\int c_{3}(144) dt =108c_{3}t-72c_{3}t+c_{4}=36c_{3}t+c_{4}.
\end{equation*}
The only restriction is 
\begin{equation*}
    W(P_{4,1},P_{4,0})(t)+12P_{4,2}(t)^2=144(36c_{3}t+c_{4})-(144t+c_{1})(36c_{3})+12c_{2}^2=0.
\end{equation*}
Hence, $c_{4}=\frac{1}{4}c_{1}c_{3}-\frac{1}{12}c_{2}^2$.

Therefore, the fundamental solution is 
\begin{equation*}
f(x,t)=x^4+c_{3}x^3+c_{2}x^2+(144t+c_{1})x+36c_{3}t+\frac{1}{4}c_{1}c_{3}-\frac{1}{12}c_{2}^2.
\end{equation*}

On the other hand, if we consider $m\geq 5$ and assume that we have a solution $f$, $\deg_{x} (f)=m$ such that $Tf=0$ then $Bf=Rf=0$ and
$R_{k}(t)=0$ for $m+1\leq k \leq 2m-1$. In particular,  $y_{k,q-1}(m)=0$ for such that $m+1\leq k \leq 2m-1$. Taking the cases $k\equiv 0$,$m+1\leq k \leq 2m-1$ we obtain $Y_{q}(m)=0$ for $q=\frac{k}{3}$. 

We have three cases:

\begin{enumerate}
     \item If $m\equiv 0$ (mod 3) then $m=3s$ with $s\geq 2$. Using Corollary \ref{leadrestric} we obtain that 
       \begin{equation*}
    V(Y_{3})=V(z_{2}+36z_{1})\cup V(z_{1})
\end{equation*}
and for $s\geq 3$
\begin{equation*}
    \bigcap_{q=s+1}^{2s-1}V(Y_{q})=\corch{V(z_{s}+36z_{s-1}) \cap \bigcap_{u=1}^{s-2}V(z_{u})}\cup \corch{\bigcap_{u=1}^{s-1}V(z_{u})}
    \cup \corch{\bigcap_{u=2}^{s}V(z_{u})},
\end{equation*}
Therefore $m\in V(z_{2}+36z_{1})\cup V(z_{1})\cup V(z_{2})$. However, by simple calculation
\begin{equation*}
 z_{1}(m)=m^2(m-1)^2, \hspace{0.1cm} z_{2}(m)=\frac{m^3(m-1)^3(m-3)(m-4)}{2}  
\end{equation*}
and 
\begin{equation*}
    z_{2}(m)+36z_{1}(m)=\frac{m^2(m-1)^2(m^4-8m^3+19m^2-12m+72)}{2}
\end{equation*}
which is a contradiction.

     \item If $m\equiv 1$ (mod 3) then $m=3s+1$ with $s\geq 2$. Using Corollary \ref{leadrestric} we obtain that 
   \begin{equation*}
    V(Y_{3})=V(z_{2}+144z_{1})\cup V(z_{1})
\end{equation*}
and for $s\geq 3$
\begin{equation*}
    \bigcap_{q=s+1}^{2s-1}V(Y_{q})=\corch{V(z_{s}+144z_{s-1}) \cap \bigcap_{u=1}^{s-2}V(z_{u})}\cup \corch{\bigcap_{u=1}^{s-1}V(z_{u})}
    \cup \corch{\bigcap_{u=2}^{s}V(z_{u})},
\end{equation*}
Therefore $m\in V(z_{2}+144z_{1})\cup V(z_{1})\cup V(z_{2})$. However, by simple calculation

\begin{equation*}
    z_{2}(m)+144z_{1}(m)=\frac{m^2(m-1)^2(m^4-8m^3+19m^2-12m+288)}{2}
\end{equation*}
which is a contradiction.
    
    \item If $m\equiv 2$ (mod 3) then $m=3s+2$ with $s\geq 1$. Using Corollary \ref{leadrestric} we obtain that 
    $m\in  \bigcap_{l=1}^{s}V(Y_{2l+1})=\bigcap_{u=1}^{s}V(z_{u})$ and therefore $m\in V(z_{1})$, which is a contradiction. \QED
\end{enumerate}
\end{demos}
\begin{coro}
There are three classes of nonzero rational solutions of the KdV-like equation \eqref{KdVlike} given by
\begin{equation*}
    u(x,t)=\frac{2}{x+c_{0}},
\end{equation*}
\begin{equation*}
    u(x,t)=\frac{2(9x^2+6c_{2}x+c_{2}^2)}{3x^3+3c_{2}x^2+c_{2}^2x+108t+3c_{0}}
\end{equation*}
and 
\begin{equation*}
    u(x,t)=\frac{2(48x^3+36c_{3}x^2+24c_{2}x+1728t+12c_{1})}{12x^4+12c_{3}x^3+12c_{2}x^2+1728tx+12c_{1}x+432c_{3}t+3c_{1}c_{3}-c_{2}^2}.
\end{equation*}
with $c_{0}, c_{1},c_{2},c_{3}$ arbitrary constants.
\end{coro}

\section*{Conclusions and Comments}

We considered a generalized bilinear equation and its associated KdV-like equation. Those types of equations are linked to the bilinear operators defined and studied in \cite{bilinearMa} and we characterized the solution based on a prime number $p=3$. We gave a positive answer to the conjecture proposed in \cite{rationalzhangMa} about two classes of rational solutions of the nonlinear KdV-like equation \eqref{KdVlike}  generated from the polynomial solutions to the generalized bilinear equation \eqref{bildiff}.

The most important fact is that the bilinear equation for $(p,k)=(2,2)$ i.e., the bilinear KdV equation $\paren{D_{2,x}D_{2,t}+D_{2,x}^4}f \cdot f=0$ has polynomial solutions of degree arbitrarily large in the variable $x$ \cite{wronskianKdV} while, we just saw, in the case $(p,k)=(3,2)$ the bilinear KdV-like equation $(D_{3,x}D_{3,t}+D_{3,x}^4)f \cdot f=0$ only has polynomial solutions of degree $1,3,4$ in the variable $x$.

The previous facts suggest that the next question in this direction is to characterize the polynomial solutions of the bilinear KdV-like equation in the case $(p,2)$ for arbitrary $p$ and to look for generalizations for several variables $k$ greater than $2$.

It is natural to think about the existence of generalized $d-$linear differential operators such as the trilinear studied in \cite{trilinear}. We think that studying their rational solutions is an interesting topic for future research. It may be possible that in the general $d-$linear case, the set of rational solutions of the KdV-like equation obtained applying the transformation \eqref{logtrans} to polynomial solutions of the generalized $d-$linear equation is finite as in the bilinear case studied in this paper for the prime number $p=3$.

\section*{Acknowledgments}
The author acknowledge the financial support provided by {\em Khalifa University} under the Grant FSU 2020-9.


\begin{thebibliography}{12}
	\providecommand{\natexlab}[1]{#1}
	\providecommand{\url}[1]{\texttt{#1}}
	\expandafter\ifx\csname urlstyle\endcsname\relax
	\providecommand{\doi}[1]{doi: #1}\else
	\providecommand{\doi}{doi: \begingroup \urlstyle{rm}\Url}\fi
	
	\bibitem[Airault et~al.(1977)Airault, McKean, and Moser]{AMM}
	H.~Airault, {H. P.} McKean, and J.~Moser.
	\newblock Rational and elliptic solutions of the korteweg-de vries equation and
	a related many-body problem.
	\newblock \emph{Communications on Pure and Applied Mathematics}, 30\penalty0
	(1):\penalty0 95--148, jan 1977.
	\newblock ISSN 0010-3640.
	\newblock \doi{10.1002/cpa.3160300106}.
	
	\bibitem[Krichever(1978)]{kricheverrat}
	I.~Krichever.
	\newblock Rational solutions of the kamdomtsev-petviashvili equation and
	integrable systems of n particles on a line.
	\newblock \emph{Functional Analysis and Its Applications - FUNCT ANAL APPL-ENGL
		TR}, 12:\penalty0 59--61, 01 1978.
	\newblock \doi{10.1007/BF01077570}.
	
	\bibitem[Zubelli(1992)]{ZUBELLI199271}
	Jorge~P Zubelli.
	\newblock Rational solutions of nonlinear evolution equations, vertex
	operators, and bispectrality.
	\newblock \emph{Journal of Differential Equations}, 97\penalty0 (1):\penalty0
	71--98, 1992.
	\newblock ISSN 0022-0396.
	\newblock \doi{https://doi.org/10.1016/0022-0396(92)90084-Z}.
	\newblock URL
	\url{https://www.sciencedirect.com/science/article/pii/002203969290084Z}.
	
	\bibitem[Sakhnovich and Zubelli(2001)]{MR1857803}
	Alexander Sakhnovich and Jorge~P. Zubelli.
	\newblock Bundle bispectrality for matrix differential equations.
	\newblock \emph{Integral Equations Operator Theory}, 41\penalty0 (4):\penalty0
	472--496, 2001.
	\newblock ISSN 0378-620X.
	\newblock \doi{10.1007/BF01202105}.
	\newblock URL \url{https://doi.org/10.1007/BF01202105}.
	
	\bibitem[Zubelli and Valerio~Silva(2000)]{Zubelli2000}
	Jorge~P. Zubelli and D.~S. Valerio~Silva.
	\newblock Rational solutions of the master symmetries of the {K}d{V} equation.
	\newblock \emph{Communications in Mathematical Physics}, 211\penalty0
	(1):\penalty0 85--109, 2000.
	\newblock ISSN 0010-3616.
	
	\bibitem[Kharif et~al.(2008)Kharif, Pelinovsky, and Slunyaev]{kharif2008rogue}
	C.~Kharif, E.~Pelinovsky, and A.~Slunyaev.
	\newblock \emph{Rogue Waves in the Ocean}.
	\newblock Advances in Geophysical and Environmental Mechanics and Mathematics.
	Springer Berlin Heidelberg, 2008.
	\newblock ISBN 9783540884194.
	\newblock URL \url{https://books.google.ae/books?id=LR-IjCcUkfIC}.
	
	\bibitem[M\"{u}ller et~al.(2005)M\"{u}ller, Garrett, and Osborne]{ocean}
	P.~M\"{u}ller, C.~Garrett, and A.~Osborne.
	\newblock Rogue waves.
	\newblock \emph{Oceanography}, September 2005.
	\newblock URL \url{https://doi.org/10.5670/oceanog.2005.30}.
	
	\bibitem[Ma(2011)]{bilinearMa}
	Wen-Xiu Ma.
	\newblock Generalized bilinear differential equations.
	\newblock \emph{Stud. Nonlinear Sci.}, 2:\penalty0 140--144, 01 2011.
	
	\bibitem[{Hirota}(2004)]{zbMATH02117215}
	Ryogo {Hirota}.
	\newblock \emph{{The direct method in soliton theory. Translated from the 1992
			Japanese original and edited by Atsushi Nagai, Jon Nimmo and Claire Gilson}},
	volume 155.
	\newblock Cambridge: Cambridge University Press, 2004.
	\newblock ISBN 0-521-83660-3; 0-511-20768-9.
	
	\bibitem[Zhang and Ma(2015)]{rationalzhangMa}
	Yi~Zhang and Wen-Xiu Ma.
	\newblock Rational solutions to a kdv-like equation.
	\newblock \emph{Applied Mathematics and Computations}, pages 252--256, 2015.
	
	\bibitem[Ma and You(2005)]{wronskianKdV}
	Wen-Xiu Ma and Yuncheng You.
	\newblock Solving the korteweg-de vries equation by its bilinear form:
	Wronskian solutions.
	\newblock \emph{Transactions of the American Mathematical Society},
	357\penalty0 (5):\penalty0 1753--1778, 2005.
	\newblock ISSN 00029947.
	\newblock URL \url{http://www.jstor.org/stable/3845133}.
	
	\bibitem[Ma(2013)]{trilinear}
	Wen-Xiu Ma.
	\newblock Trilinear equations, bell polynomials, and resonant solutions.
	\newblock \emph{Frontiers of Mathematics in China}, pages 1673--3576, 10 2013.
	
\end{thebibliography}
\end{document}